\documentclass[12pt]{article}
\usepackage{amssymb,amsmath}
\usepackage{cases}
\usepackage{amsfonts}
\usepackage{color,xcolor}
\usepackage[left=2.0cm,right=2.0cm,top=2.0cm,bottom=2.0cm]{geometry}
\usepackage[colorlinks,citecolor=blue,urlcolor=blue]{hyperref}

\newtheorem{theorem}{Theorem}[section]

\newtheorem{lemma}{Lemma}[section]

\newtheorem{definition}{Definition}[section]
\newtheorem{remark}{Remark}[section]

\newcommand{\bal}{\begin{align}}
\newcommand{\bbal}{\begin{align*}}
\newcommand{\beq}{\begin{equation}}
\newcommand{\eeq}{\end{equation}}
\newcommand{\bca}{\begin{cases}}
\newcommand{\eca}{\end{cases}}
\def\div{\mathord{{\rm div}\ }}
\newcommand{\pa}{\partial}
\newcommand{\fr}{\frac}
\newcommand{\na}{\nabla}
\newcommand{\De}{\Delta}

\newcommand{\cd}{\cdot}
\newcommand{\ep}{\varepsilon}
\newcommand{\dd}{\ \mathrm{d}}
\newcommand{\B}{\dot{B}}

\newcommand{\R}{\mathbb{R}}
\newcommand{\les}{\lesssim}

\newcommand{\bi}{\Big}
\newcommand{\uu}{\mathbf{u}}
\newcommand{\vv}{\mathbf{v}}
\newcommand{\w}{\mathbf{w}}
\newcommand{\ccc}{\mathbf{c}}
\newcommand{\U}{\mathbf{U}}
\newcommand{\W}{\mathbf{W}}
\newcommand{\F}{\mathbf{F}}
\newcommand{\G}{\mathbf{G}}
\newcommand{\HH}{\mathbf{H}}

\begin{document}
\title{A class large solution of 3D  incompressible micropolar fluid system}

\author{Jinlu Li$^{1}$\footnote{E-mail: lijinlu@gnnu.edu.cn}, \qquad Weipeng Zhu$^{2}$\footnote{E-mail: mathzwp2010@163.com(Corresponding author)}\\
\small $^1$\it School of Mathematics and Computer Sciences, Gannan Normal University, Ganzhou 341000, China\\
\small $^2$\it School of Mathematics and Information Science, Guangzhou University, Guangzhou 510006, China}

\date{}

\maketitle\noindent{\hrulefill}

{\bf Abstract:} In this paper, we construct a class global large solution to the 3D  incompressible micropolar fluid system. Precisely speaking, by choosing a special initial data which can be arbitrarily large in $\dot{B}^{-1}_{\infty,\infty}$, the system has a unique global solution.


\vskip0mm\noindent{\hrulefill}

\section{Introduction}\label{sec1}
This paper focuses on the following 3D incompressible micropolar fluid system given by
\begin{eqnarray}\label{3d-m}
        \left\{\begin{array}{ll}
          \partial_t\uu+\uu\cd\na \uu-(\chi+\nu)\De \uu+\na p-2\chi\na\times \w=0,& x\in \R^3,t>0,\\
          \partial_t\w+\uu\cd\na \w-\mu\Delta \w-\kappa\nabla\operatorname{div}\w+4\chi \w-2\chi\na\times \uu=0,& x\in \R^3,t>0,\\
          \div \uu=0,& x\in \R^3,t\geq0,\\
                   (\uu,\w)|_{t=0}=(\uu_0,\w_0),& x\in \R^3.\end{array}\right.
        \end{eqnarray}
Here $\uu=\uu(x,t)$, $\w=\w(x,t)$ and $p=p(x,t)$ are unknown functions representing the linear velocity field, the micro-rotation velocity field and the pressure field of the fluid, respectively, and $\mu$, $\kappa$, $\nu$ and $\chi$ are nonnegative constants reflecting various viscosity of the fluid.

  The model \eqref{3d-m} was first proposed by Eringen in the pioneering work \cite{E66}. It is an essential
  modification to the classical Navier-Stokes equations for the purpose to better describe the motion of various real fluids
  consisting of rigid but randomly oriented particles (e.g., blood) by considering the effect of micro-rotation of the particles
  suspended in the fluid. A fluid possessing such a property is called a {\em micropolar fluid}, so that the model \eqref{3d-m} is referred
  to as {\em micropolar fluid system} in the literature.
  Actually, there are several experiments indicate that solutions of the micropolar fluid system do better mimic behavior of such kind of fluids than solutions of the Navier--Stokes equations, cf., \cite{MP84, P69, PRU74} and references therein. We refer the reader to see the references \cite{L99, L03} for more physical background of the above model.

The micropolar fluid system has recently drawn much attention and many interesting results have been obtained in the literature. For instance, the first result on existence and uniqueness of solutions of the problem \eqref{3d-m} was obtained by Galdi and Rionero in the reference
  \cite{GR77}. Existence of global weak solutions of the problem \eqref{3d-m} was established by Lukaszewicz \cite{Luk89} and Boldrini and
  Rojas-Medar \cite{RojB98}. We refer the reader to see \cite{BolRF03,Roj97,RojO05} for existence and uniqueness of strong solutions to the system \eqref{3d-m} and more complex systems such
  as the magneto-micropolar fluid system. Well-posedness of the micropolar fluid system
  in various function spaces has also been well studied by many authors. For example,
  Ferreira and Villamizar-Roa \cite{FV07} proved well-posedness of a generalized incompressible micropolar fluid system in pseudo-measure spaces.
  In \cite{CM12}, Chen and Miao established global well-posedness of the system \eqref{3d-m} for small initial data in the
  Besov spaces $\dot{B}_{p,\infty}^{-1+\frac{3}{p}}(\mathbb{R}^3)$ for $p\in [1,6)$.
  In \cite{ZhuZ18}, Zhu and Zhao proved that the micropolar fluid system is globally well-posed in the Fourier-Besov spaces $\dot{FB}_{p,r}^{2-\frac{3}{p}}(\mathbb{R}^3)$ for
  $p\in (1,\infty]$ and $r\in [1,\infty)$ with small initial data. Recently, the corresponding author of the present paper \cite{Zhu19} showed that this problem is well-posed in $\dot{FB}_{1,r}^{-1}(\mathbb{R}^3)$ for
  $1\leqslant r\leqslant 2$, while ill-posed for $2<r\leqslant\infty$. We also refer the reader to see the references \cite{DonLW17,VilR08,WanW17,Yam05}
  for other related work.

  As mentioned above, the global regularity for the micropolar fluid system \eqref{3d-m} with
large initial data is still a challenging open problem. Our main goal is to construct a class global large solution to problem \eqref{3d-m}. For simplicity, throughout this paper we only consider the situation with $\mu=\kappa=1$ and $\chi=\nu={1}/{2}$. The main result is stated as follows.

\begin{theorem}\label{the1.1} Let $p\in(4,6)$. Assume that the initial data fulfills ${\rm{div}}\vv_0=0$ and
$$\uu_0=\U_0+\vv_0\quad \mbox{and}\quad \w_0=\W_0+\ccc_0$$
where
\begin{eqnarray*}
&\U_0=
\begin{pmatrix}
\pa_2a_0 \\ -\pa_1a_0\\ 0
\end{pmatrix}
\quad\mbox{and}\quad
\W_0=
\begin{pmatrix}
0\\ 0 \\ a_0
\end{pmatrix}
\end{eqnarray*}
with
\begin{eqnarray}\label{con-Equ1.2}
\mathrm{supp} \ \hat{a}_0(\xi)\subset\mathcal{C}:=\Big\{\xi \big| \ |\xi_1+\xi_2|\leq \ep, \ 1\leq |\xi_h|\leq 2, \ \ep\leq |\xi_3|\leq 2\ep\Big\} .
\end{eqnarray}
There exists a sufficiently small positive constant $\delta$, and a universal constant $C$ such that if
\begin{align}\label{condition}
(||\vv_0,\ccc_0||_{\B^{-1+\frac3p}_{p,1}}+\ep||\hat{a}_0||^2_{L^{\frac{p}{p-1}}}+\ep||\hat{a}_0||_{L^{\frac{p}{p-1}}} )\exp\bi\{C\bi( ||\hat{a}_0||^2_{L^1}+||\hat{a}_0||_{L^1}\bi)\bi\}\leq \delta,
\end{align}
then the system \eqref{3d-m} has a unique global solution.
\end{theorem}

\begin{remark}\label{rem1.1}
Let $\vv_0=\ccc_0=0$ and $a_0(x_1,x_2,x_3)=\ep^{-2}(\log\log\frac1\ep)^{\frac12}  \chi(x_1,x_2)\phi(x_3)$.
where the smooth functions $\chi,\phi$ satisfying
\begin{align*}
\mathrm{supp} \hat{\chi}\in \mathcal{\widetilde{C}},\quad \hat{\chi}(\xi)\in[0,1]\quad\mbox{and} \quad \hat{\chi}(\xi)=1 \quad\mbox{for} \quad \xi\in\mathcal{\widetilde{C}}_1,
\end{align*}
and
\begin{align*}
\mathrm{supp} \hat{\phi}(\xi') \in [\ep,2\ep],\quad \hat{\phi}(\xi') \in [0,1] \quad\mbox{and} \quad \hat{\phi}(\xi)=1 \quad\mbox{for}  \quad \xi'\in[\frac54\ep,\frac74\ep],
\end{align*}
where
\begin{align*}
&\mathcal{\widetilde{C}}\triangleq\Big\{\xi\in\R^2: \ |\xi_1+\xi_2|\leq \ep,\ 1\leq\xi^2_1+\xi^2_2\leq 2\Big\},
\\&\mathcal{\widetilde{C}}_1\triangleq\Big\{\xi\in\R^2: \ |\xi_1+\xi_2|\leq \frac12\ep,\ \frac54\leq\xi^2_1+\xi^2_2\leq \frac74\Big\} .
\end{align*}
Then, direct calculations show that the left side of \eqref{condition} becomes
\begin{align*}
C\ep^{\frac{p-4}{p}} \log\log \frac1\ep e^{C\log\log\frac1\ep}.
\end{align*}
Therefore, choosing $\ep$ small enough, we deduce that the system \eqref{3d-m} has a global solution.

Note that $\omega_0=\partial_2 \uu_0^1-\partial_1\uu_0^2=(\partial_1^2+\partial_2^2)a_0$ and $\hat{a}_0\geq 0$,  we can deduce that $\hat{\omega}_0=-(\xi_1^2+\xi_2^2)\hat{a}_0$ which implies $\|\omega_0\|_{L^\infty}=\|\hat{\omega}_0\|_{L^1}$. Moreover, we also have $\|\hat{\omega}_0\|_{L^1}\gtrsim (\log\log \frac1\ep)^\frac12$. Using the inequality $\|\omega_0\|_{L^\infty}\lesssim \|\uu_0\|_{L^\infty}$, then we have $||\uu_0||_{L^\infty}\gtrsim (\log\log \frac1\ep)^\frac12$. Therefore, we obtain $\|\uu_0\|_{\dot{B}^{-1}_{\infty,\infty}}\approx||\uu_0||_{L^\infty}\gtrsim (\log\log \frac1\ep)^\frac12$. Obviously, we also have $\|\w_0\|_{\dot{B}^{-1}_{\infty,\infty}}\approx||\w_0||_{L^\infty}\gtrsim (\log\log \frac1\ep)^\frac12$.
\end{remark}

\begin{remark}\label{rem1.1}
The result in Theorem \ref{the1.1} is also valid when $\kappa=0$ ($\mu,\chi,\nu>0$).
\end{remark}

\section{Littlewood-Paley Analysis}
Throughout this paper, we will denote by $C$ any constant which may change from line to line and write $A\lesssim B$ if $A\leq CB$. $A\approx B$ means that $A\lesssim B$ and $B\lesssim A$. We also shall use the abbreviated notation $||f_1,\cdots,f_n||_{X}=||f_1||_{X}+\cdots+||f_n||_{X}$ for some Banach space $X$.

Next, we recall the Littlewood-Paley theory, the definition of homogeneous Besov spaces and some useful properties.

Let us start by introducing the Littlewood-Paley decomposition. Choose a radial function $\varphi\in \mathcal{S}(\mathbb{R}^d)$ supported in ${\mathcal{C}}=\{\xi\in\mathbb{R}^d:\frac34\leq |\xi|\leq \frac83\}$ such that
\begin{align*}
\sum_{j\in \mathbb{Z}}\varphi(2^{-j}\xi)=1 \quad \mathrm{for} \ \mathrm{all} \ \xi\neq0.
\end{align*}
The frequency localization operator $\dot{\Delta}_j$ and $\dot{S}_j$ are defined by
\begin{align*}
\dot{\Delta}_jf=\varphi(2^{-j}D)f=\mathcal{F}^{-1}(\varphi(2^{-j}\cdot)\mathcal{F}f) \quad \mbox{and}\quad\dot{S}_jf=\sum_{k\leq j-1}\dot{\Delta}_kf \quad \mathrm{for} \quad j\in\mathbb{Z}.
\end{align*}
With a suitable choice of $\varphi$, one can easily verify that
\begin{align*}
\dot{\Delta}_j\dot{\Delta}_kf=0 \quad \mathrm{if} \quad |j-k|\geq2\quad \mbox{and}\quad\dot{\Delta}_j(\dot{S}_{k-1}f\dot{\Delta}_kf)=0 \quad  \mathrm{if} \quad  |j-k|\geq5.
\end{align*}
Next we recall Bony's decomposition from \cite{Bahouri2011}:
\begin{align*}
uv=\dot{T}_uv+\dot{T}_vu+\dot{R}(u,v),
\end{align*}
with
\begin{align*}
\dot{T}_uv=\sum_{j\in\mathbb{Z}}\dot{S}_{j-1}u\dot{\Delta}_jv, \quad \quad \dot{R}(u,v)=\sum_{j\in\mathbb{Z}}\dot{\Delta}_ju\widetilde{\Delta}_jv, \quad \quad \widetilde{\Delta}_jv=\sum_{|j'-j|\leq 1}\dot{\Delta}_{j'}v.
\end{align*}
\begin{definition}
We denote by $\mathcal{Z}'(\mathbb{R}^d)$ the dual space of $\mathcal{Z}(\mathbb{R}^d)$, where we set $$\mathcal{Z}(\mathbb{R}^d)=\bi\{f\in \mathcal{S}(\mathbb{R}^d):D^{\alpha}\hat{f}(0)=0;\ \forall \alpha\in \mathbb{N}^d\bi\}.$$
\end{definition}
Then we have the formal homogenous Littlewood-Paley decomposition $$f=\sum\limits_{j\in \mathbb{Z}}\dot{\Delta}_jf,\quad \forall f\in\mathcal{Z}'(\mathbb{R}^d).$$

The operators $\dot{\Delta}_j$ help us recall the definition of the homogenous Besov space and the related lemma (see \cite{Bahouri2011}).

\begin{definition}
Let $s\in \mathbb{R}$ and $1\leq p,r\leq \infty$. The homogeneous Besov space $\dot{B}^s_{p,r}$ is defined by
\begin{align*}
\dot{B}^s_{p,r}=\bi\{f\in \mathcal{Z}'(\mathbb{R}^d):||f||_{\dot{B}^s_{p,r}}\triangleq \Big|\Big|(2^{ks}||\dot{\Delta}_k f||_{L^p})_{k\in\mathbb{Z}}\Big|\Big|_{\ell^r}<+\infty\bi\}.
\end{align*}
\end{definition}
It should be noted that a distribution $f\in\dot{B}^s_{p,r}$ if and only if there exist a constant $C>0$ and a non-negative sequence $\{d_{k}\}_{k\in\mathbb{Z}}$ such that
\begin{eqnarray*}
\forall k\in\mathbb{Z},\quad \|\dot{\Delta}_kf\|_{L^p}\leq Cd_{k} 2^{-ks}\|f\|_{{\dot{B}_{p,r}^{s}}} \quad \mbox{with} \quad \|d_{k}\|_{\ell^r}=1.
\end{eqnarray*}

\begin{lemma}\label{lem2.1}
Let $t>0$ and $1\leq p,p_1,p_2,r,r_1,r_2\leq \infty$ satisfying
\bbal
\frac1p=\frac{1}{p_1}+\frac{1}{p_2},\qquad \frac1r=\frac{1}{r_1}+\frac{1}{r_2}, \qquad t=t_1+t_2.
\end{align*}
Then there exists a positive constant $C$ such that
\begin{align}
&||\dot{T}_fg||_{\B^s_{p,r}}\leq C||f||_{L^{p_1}}||g||_{\B^s_{p_2,r}},\label{01}\\
&||\dot{T}_fg||_{\B^{s-t}_{p,r}}\leq C||f||_{\B^{-t}_{p_1,r_1}}||g||_{\B^s_{p_2,r_2}},\label{02}
\\&||\dot{R}(f,g)||_{\B^t_{p,r}}\leq C||f||_{B^{t_1}_{p_1,r_1}}||g||_{B^{t_2}_{p_2,r_2}}.\label{03}
\end{align}
\end{lemma}

Next, we present the following product estimate which will be used in the sequel.
\begin{lemma}\cite{Danchin2001}\label{lem2.2}
Let $2\leq p\leq \infty$, $s_1\leq \frac 3p$ and $s_2\leq \frac 3p$ with $s_1+s_2>3\max\{0,\frac2p-1\}$. Then there holds
\begin{align*}
||fg||_{\dot{B}^{s_1+s_2-\frac 3p}_{p,1}}&\leq C||f||_{\dot{B}^{s_1}_{p,1}}||g||_{\dot{B}^{s_2}_{p,1}}.
\end{align*}
\end{lemma}

Finally, we recall the optimal regularity estimates for the heat equations.
\begin{lemma}\cite{Danchin2005}\label{lem2.3}
Let $s\in \mathbb{R}$, $1\leq p,r\leq \infty$  and $1\leq q_1\leq q_2\leq \infty$. Assume that $u_0\in \dot{B}^s_{p,r}$ and $G\in {\widetilde{L}}^{q_1}_T(\dot{B}^{s+\frac{2}{q_1}-2}_{p,r})$. Then the heat equations
\begin{align*}
\left\{\begin{array}{ll}
\partial_tu-\Delta u=G,\\
 u(0,x)=u_0,
\end{array}\right.
\end{align*}
has a unique solution $u\in \widetilde{L}^{q_2}_T(\dot{B}^{s+\fr{2}{q_2}}_{p,r}))$ satisfying for all $T>0$
\begin{align*}
||u||_{\widetilde{L}^{q_2}_T(\dot{B}^{s+\fr{2}{q_2}}_{p,r})}\les||u_0||_{\dot{B}^s_{p,r}}+||G||_{{\widetilde{L}}^{q_1}_T(\dot{B}^{s+\frac{2}{q_1}-2}_{p,r})}.
\end{align*}
\end{lemma}

\section{Proof of the main theorem}

In this section, we will give the proof of Theorem \ref{the1.1}.

Let $(a,m)$ be the solutions of the following system
\begin{eqnarray}\label{m-l1}
        \left\{\begin{array}{ll}
          \pa_ta-\De a- m=0,\\
          \pa_tm-\De m+2 m+ \De a=0,\\
          (a,m)|_{t=0}=(a_0,a_{0}),\end{array}\right.
        \end{eqnarray}
and the action of its Green matrix which is denoted by  $G(t,x)$. We can show that $G(t,x)$ has similar properties in common with the heat kernel, that is
$$|\hat{G}(t,\xi)|\leq e^{-c|\xi|^2t},\quad |\xi|>0.$$
In fact, we have

\bal\label{greem-1}
\widehat{Gf}(t,\xi)=e^{-A(\xi)t}\hat{f}(\xi), \quad f=(a_0,a_0)^T.
\end{align}
where
\begin{eqnarray*}
A(\xi)=
\begin{pmatrix}
|\xi|^2& -1\\ -|\xi|^2& |\xi|^2+2
\end{pmatrix}
.
\end{eqnarray*}
Then for $|\xi|>0$,
\begin{align*}
|\widehat{Gf}(t,\xi)| = |e^{-A(\xi)t}\hat{f}(\xi)| \leq e^{-(|\xi|^2+1-\sqrt{|\xi|^2+1})t} |\hat{f}(\xi)|\leq e^{-c|\xi|^2t}|\hat{f}(\xi)|.
\end{align*}
If the initial data $a_0$ satisfies some assumptions, we can find small quantities along with $a_0$ which is the key to construct global solutions in our paper. That is, for $p\geq 2$
\bbal
||(\pa_1+\pa_2)a||_{L^p}+||(\pa_1+\pa_2)m||_{L^p}+||\pa_3a||_{L^p}+||\pa_3m||_{L^p}\lesssim \ep e^{-ct}||\hat{a}_0||_{L^{\frac{p}{p-1}}}.
\end{align*}
Setting
\bbal
&\U=
\begin{pmatrix}
\pa_2a \\ -\pa_1a\\ 0
\end{pmatrix}
 \quad\mbox{and}\quad
\W=
\begin{pmatrix}
0\\ 0 \\ m
\end{pmatrix},
\end{align*}
we can deduce from \eqref{m-l1} that
\begin{eqnarray}\label{m-l2}
        \left\{\begin{array}{ll}
          \pa_t\U-\De \U- \na\times \W=0,\\
          \pa_t\W-\De \W-\nabla\operatorname{div}\W+2 \W- \na\times \U=\F,\\
          \div \U=0, \quad (\U,\W)|_{t=0}=(\U_0,\W_{0}),\end{array}\right.
        \end{eqnarray}
where
\bbal
\F=
\begin{pmatrix}
-\pa_1\pa_3(a+m) \\ -\pa_2\pa_3(a+m)\\ -\pa_3^2 (a+m)
\end{pmatrix}
 .
\end{align*}
Denoting $\vv=\uu-\U$ and $\ccc=\w-\W$, the system \eqref{3d-m} can be written as follows
\begin{eqnarray}\label{c-3d-m}
        \left\{\begin{array}{ll}
          \partial_t\vv+\vv\cd\na \vv-\De \vv+\na p-\na\times \ccc=\G-\U\cd\na \vv-\vv\cd\na \U,\\
          \partial_t\ccc+\vv\cd\na \ccc-\Delta \ccc-\nabla\operatorname{div}\ccc+2 \ccc-\na\times \vv=\HH-\F-\U\cd\na \ccc-\vv\cd\na \W,\\
          \div \vv=0,\\
                   (\vv,\ccc)|_{t=0}=(\vv_0,\ccc_0).
          \end{array}\right.
        \end{eqnarray}
where
\bbal
&\G=-\U\cd\na \U,\qquad  \HH=-\U\cd\na \W.
\end{align*}

To element the pressure term, applying the Leray operator $\mathbb{P}$ to the first equation of \eqref{c-3d-m}, one has
\begin{eqnarray}\label{3.1}
        \left\{\begin{array}{ll}
          \partial_t\vv-\De \vv-\na\times \ccc=\mathbb{P}(-\vv\cd\na \vv+\G-\U\cd\na \vv-\vv\cd\na \U),\\
          \partial_t\ccc-\Delta \ccc-\nabla\operatorname{div}\ccc+2 \ccc-\na\times \vv=-\vv\cd\na \ccc+\HH-\F-\U\cd\na \ccc-\vv\cd\na \W,\\
          \div \vv=0,\\
                   (\vv,\ccc)|_{t=0}=(\vv_0,\ccc_0).
          \end{array}\right.
\end{eqnarray}
Note that it has been shown in \cite{CM12} (Proposition 3.5) that the Green matrix of the linear system of \eqref{3.1} has similar properties in common with the heat kernel by making a suitable transformation of the solution. Then, invoking Lemma \ref{lem2.3} to the above system \eqref{3.1} yields
\bal\label{3.2}
&||\vv,\ccc||_{L^\infty_t(\B^{-1+\frac3p}_{p,1})}+||\vv,\ccc||_{L^1_t(\B^{1+\frac3p}_{p,1})}
\\\les&~ ||\vv_0,\ccc_0||_{\B^{-1+\frac3p}_{p,1}}+\underbrace{\int^t_0||\vv\cd\na \vv||_{\B^{\frac3p-1}_{p,1}}+||\vv\cd\na \ccc||_{\B^{\frac3p-1}_{p,1}}\dd \tau}_{I_1}+\underbrace{\int^t_0||\G,\HH,\F||_{\B^{\frac3p-1}_{p,1}}\dd \tau}_{I_2}\nonumber\\
\quad&+\underbrace{\int^t_0||\U\cd\na \vv+\vv\cd\na \U||_{\B^{\frac3p-1}_{p,1}}+||\U\cd\na \ccc+\vv\cd\na \W||_{\B^{\frac3p-1}_{p,1}}\dd \tau}_{I_3},
\end{align}
where we have used the fact that $\mathbb{P}$ is a smooth homogeneous of degree 0 Fourier multipliers which maps $\dot{B}_{p,1}^{\frac3p-1}$ to itself.

For the term $I_1$, using the product estimate (see Lemma \ref{lem2.2}), we obtain
\bal\label{3.3}
I_1 & \les~\int_0^t(||\vv||_{\B^{-1+\frac3p}_{p,1}}+||\ccc||_{\B^{-1+\frac3p}_{p,1}})(||\vv||_{\B^{1+\frac3p}_{p,1}}+|\ccc||_{\B^{1+\frac3p}_{p,1}})\dd \tau \nonumber\\
&\les~||\vv,\ccc||_{L^\infty_t(\B^{-1+\frac3p}_{p,1})}||\vv,\ccc||_{L^1_t(\B^{1+\frac3p}_{p,1})}.
\end{align}

For the term $I_2$, we estimate $\G, \HH, \F$ one by one.  Direct calculations show that $\G^3=\HH^1=\HH^2=0$ and
\bbal
\G^1&=-\U\cd\na \U^1=\pa_1a\pa_2\pa_2a-\pa_2a\pa_1\pa_2a
\\&=(\pa_1+\pa_2)a\pa_2\pa_2a-\pa_2a\pa_2(\pa_1+\pa_2)a,
\end{align*}
\bbal
\G^2&=-\U\cd\na \U^2=-\pa_1a\pa_2\pa_1a+\pa_2a\pa_1\pa_1a
\\&=-(\pa_1+\pa_2)a\pa_1\pa_2a+\pa_2a\pa_1(\pa_1+\pa_2)a,
\end{align*}
\bbal
\HH^3&=-\U\cd\na \W^3=\pa_1a\pa_2m-\pa_2a\pa_1m
\\&=(\pa_1+\pa_2)a\pa_2m-\pa_2a(\pa_1+\pa_2)m.
\end{align*}
Using the Lemma \ref{lem2.2} again gives
\bal\label{3.4}
\int^t_0||\G||_{\B^{-1+\frac3p}_{p,1}}\dd \tau\les&~ \int^t_0||(\pa_1+\pa_2)a||_{\B^{\frac3p}_{p,1}}||a||_{\B^{1+\frac3p}_{p,1}}\dd\tau
\les~ \ep||\hat{a}_0||^2_{L^{\frac{p}{p-1}}},
\end{align}
and
\bal\label{3.4-1}
\int^t_0||\HH||_{\B^{-1+\frac3p}_{p,1}}\dd \tau\les&~ \int^t_0||(\pa_1+\pa_2)a||_{\B^{\frac3p}_{p,1}}||m||_{\B^{\frac3p}_{p,1}}\dd\tau
\les~ \ep||\hat{a}_0||^2_{L^{\frac{p}{p-1}}},
\end{align}
The term $\F$ can be calculated as follows
\bal\label{3.4-2}
\int^t_0||\F||_{\B^{-1+\frac3p}_{p,1}}\dd \tau\les&~ \int^t_0||\pa_3(a+m)||_{\B^{\frac3p}_{p,1}}\dd\tau
\les~ \ep||\hat{a}_0||_{L^{\frac{p}{p-1}}}.
\end{align}
Putting \eqref{3.4}--\eqref{3.4-2} together yields
\bbal
I_2\lesssim \ep||\hat{a}_0||^2_{L^{\frac{p}{p-1}}}+\ep||\hat{a}_0||_{L^{\frac{p}{p-1}}}.
\end{align*}

To deal with the term $I_3$, by Bony's decomposition, one has
\begin{equation*}
\vv\cd\na \U^i=\sum^{3}_{j=1}\pa_j(\U^i\vv^j)=\sum^{3}_{j=1}[\pa_j(\dot{T}_{\U^i} \vv^j)+\dot{T}_{\vv^j} \pa_j\U^i+\pa_j\dot{R}(\U^i,\vv^j)]
\end{equation*}
and
\begin{equation*}
\U\cd\na \vv^i=\sum^{3}_{j=1}\pa_j(\U^j\vv^i)=\sum^{3}_{j=1}[\dot{T}_{\U^j} \pa_j\vv^i+\pa_j(\dot{T}_{\vv^i} \U^j)+\pa_j\dot{R}(\U^j,\vv^i)].
\end{equation*}
Using \eqref{01}--\eqref{03} from Lemma \ref{lem2.1}, respectively, we obtain
\bal
&||\pa_j(\dot{T}_{\U^i} \vv^j),\dot{T}_{\U^j} \pa_j\vv^i||_{\B^{-1+\frac3p}_{p,1}}\les ~||\U||_{L^\infty}||\vv||_{\B^{\frac3p}_{p,1}},\label{y1}\\
&||\dot{T}_{\vv^j} \pa_j\U^i,\pa_j(\dot{T}_{\vv^i} \U^j)||_{\B^{-1+\frac3p}_{p,1}}\les~||\U||_{\B^1_{\infty,\infty}}||\vv||_{\B^{-1+\frac3p}_{p,1}},\label{y2}\\
&||\dot{R}(\U^i,\vv^j),\dot{R}(\U^j,\vv^i)||_{\B^{\frac3p}_{p,1}}\les~||\U||_{\B^1_{\infty,\infty}}||\vv||_{\B^{-1+\frac3p}_{p,1}}.\label{y3}
\end{align}
Combining \eqref{y1}--\eqref{y3} implies
\bal\label{3.5}
&\int^t_0||\vv\cd\na \U+\U\cd\na \vv||_{\B^{-1+\frac3p}_{p,1}}\dd \tau \nonumber \\\les&~ \int^t_0||\vv||_{\B^{\frac3p}_{p,1}}||\U||_{L^\infty}\dd \tau+\int^t_0||\vv||_{\B^{-1+\frac3p}_{p,1}}||\U||_{\B^1_{\infty,\infty}}\dd \tau\nonumber\\
\les&~ \int^t_0||\vv||_{\B^{-1+\frac3p}_{p,1}}||\U||_{\B^1_{\infty,\infty}}\dd \tau+\int^t_0||\U||_{L^\infty}||\vv||^{\frac{1}{2}}_{\B^{-1+\frac3p}_{p,1}}||\vv||^{\frac{1}{2}}_{\B^{1+\frac3p}_{p,1}}\dd \tau\nonumber\\
\leq&~ C\int^t_0\bi(||\U||^2_{L^\infty}+||\U||_{\B^{1}_{\infty,\infty}}\bi)||\vv||_{\B^{-1+\frac3p}_{p,1}}\dd \tau+\frac12||\vv||_{L^1_t(\B^{1+\frac3p}_{p,1})}.
\end{align}
Similar argument as \eqref{3.5}, we also have
\bal\label{3.5-1-1}
&\int^t_0||\vv\cd\na \W+\U\cd\na \ccc||_{\B^{-1+\frac3p}_{p,1}}\dd \tau \nonumber
\\\leq&~ C\int^t_0\bi(||\U,\W||^2_{L^\infty}+||\U,\W||_{\B^{1}_{\infty,\infty}}\bi)||\vv,\ccc||_{\B^{-1+\frac3p}_{p,1}}\dd \tau+\frac12||\vv,\ccc||_{L^1_t(\B^{1+\frac3p}_{p,1})}.
\end{align}
Putting the estimates \eqref{3.3}, \eqref{3.4} and \eqref{3.5} together with \eqref{3.2} yields
\bal\label{yyh}
||\vv,\ccc&||_{L^\infty_t(\B^{-1+\frac3p}_{p,1})}+||\vv,\ccc||_{L^1_t(\B^{1+\frac3p}_{p,1})}
\les~||\vv_0,\ccc_0||_{\B^{-1+\frac3p}_{p,1}}+||\vv,\ccc||_{L^\infty_t(\B^{-1+\frac3p}_{p,1})}||\vv,\ccc||_{L^1_t(\B^{1+\frac3p}_{p,1})}
\nonumber\\
\quad& +\int^t_0\bi(||\U,\W||^2_{L^\infty}+||\U,\W||_{\B^{1}_{\infty,\infty}}\bi)||\vv,\ccc||_{\B^{-1+\frac3p}_{p,1}}\dd \tau+
\ep||\hat{a}_0||^2_{L^{\frac{p}{p-1}}}+\ep||\hat{a}_0||_{L^{\frac{p}{p-1}}}.
\end{align}
Now, we define
\bbal
\Gamma\triangleq\sup\bi\{t\in[0,T^*):||\vv,\ccc||_{L^\infty_t(\B^{-1+\frac3p}_{p,1})}\leq \eta\ll1\bi\},
\end{align*}
where $\eta$ is a small enough positive constant which will be determined later on.

Thus, for all $t\in[0,\Gamma]$, choosing $\eta$ small enough, we infer from \eqref{yyh} that
\bal\label{es-con1}
||\vv,\ccc||_{L^\infty_t(\B^{-1+\frac3p}_{p,1})}+||\vv,\ccc||_{L^1_t(\B^{1+\frac3p}_{p,1})}\les&~ ||\vv_0,\ccc_0||_{\B^{-1+\frac3p}_{p,1}}+\ep||\hat{a}_0||^2_{L^{\frac{p}{p-1}}}+\ep||\hat{a}_0||_{L^{\frac{p}{p-1}}}
\nonumber\\
\quad& +\int^t_0\bi(||\U,\W||^2_{L^\infty}+||\U,\W||_{\B^{1}_{\infty,\infty}}\bi)||\vv,\ccc||_{\B^{-1+\frac3p}_{p,1}}\dd \tau.
\end{align}
Notice that
$$\int^t_0\bi(||\U,\W||^2_{L^\infty}+||\U,\W||_{\B^{1}_{\infty,\infty}}\bi)\dd \tau\les ||\hat{a}_0||^2_{L^1}+||\hat{a}_0||_{L^1},$$
by Gronwall's inequality, we have for all $t\in[0,\Gamma]$
\bal\label{es-con2}
&||\vv,\ccc||_{L^\infty_t(\B^{\frac3p-1}_{p,1})}+||\vv,\ccc||_{L^1_t(\B^{\frac3p+1}_{p,1})}\nonumber
\\\leq&~ C(||\vv_0,\ccc_0||_{\B^{-1+\frac3p}_{p,1}}+\ep||\hat{a}_0||^2_{L^{\frac{p}{p-1}}}+\ep||\hat{a}_0||_{L^{\frac{p}{p-1}}} )\exp\bi\{C\bi( ||\hat{a}_0||^2_{L^1}+||\hat{a}_0||_{L^1}\bi)\bi\}\nonumber\\
\leq&~C\delta
\end{align}
provided that the condition \eqref{condition} holds.

Choosing $\eta=2C\delta$, we can get
\bbal
||v||_{L^\infty_t(\B^{-1+\frac3p}_{p,1})}&\leq \fr\eta2 \quad\mbox{for}\quad t\leq \Gamma.
\end{align*}

So if $\Gamma<T^*$, due to the continuity of the solutions, we can obtain that there exists $0<\epsilon\ll1$ such that
\bbal
||\vv,\ccc||_{L^\infty_t(\B^{-1+\frac3p}_{p,1})}&\leq \eta \quad\mbox{for}\quad t\leq \Gamma+\epsilon<T^*,
\end{align*}
which is contradiction with the definition of $\Gamma$.

Thus, we can conclude $\Gamma=T^*$ and
\bbal
||\vv,\ccc||_{L^\infty_t(\B^{-1+\frac3p}_{p,1})}&\leq C<\infty \quad\mbox{for all}\quad t\in(0,T^*),
\end{align*}
which implies that $T^*=+\infty$.

\section*{Acknowledgments} J. Li is supported by the National Natural Science Foundation of China (Grant No.11801090). W. Zhu is partially supported by the National Natural Science Foundation of China (Grant No.11901092) and Natural Science Foundation of Guangdong Province (No.2017A030310634).

\end{document}